\documentclass[11pt, reqno]{amsart}

\usepackage{esint}
\usepackage{amstext}
\usepackage{amsthm}
\usepackage{amsmath}
\usepackage{amssymb}
\usepackage{latexsym}
\usepackage{amsfonts}
\usepackage{graphicx}
\usepackage{color}

\usepackage[%pagebackref,
hypertexnames=false, colorlinks, citecolor=black, linkcolor=blue, urlcolor=red]{hyperref}

\newcommand{\ci}[1]{_{ {}_{\scriptstyle #1}}}

\newcommand{\norm}[1]{\ensuremath{\left\|#1\right\|}}

\newenvironment{entry}
{\begin{list}{X}%
  {%
      \setlength{\labelwidth}{55pt}%
      \setlength{\leftmargin}{\labelwidth}%\labelsep}%
      \addtolength{\leftmargin}{\labelsep}%
      \setlength{\itemsep}{.4pc}
   }%
}%
{\end{list}}

\newcommand{\R}{\mathbb{R}}

\newcommand{\C}{\mathbb{C}}

\newcommand{\F}{\mathbb{F}}
\newcommand{\E}{\mathbb{E}}

\newcommand{\bI}{\mathbf{I}}

\newcommand{\dd}{\mathrm{d}}

\newcommand{\cF}{\mathcal{F}}
\newcommand{\cX}{\mathcal{X}}

\newcommand{\cA}{\mathcal{A}}
\newcommand{\cJ}{\mathcal{J}}

\newcommand{\cD}{\mathcal{D}}

\newcommand{\1}{\mathbf{1}}

\newcommand{\fdot}{\,\cdot\,}
\newcommand{\wt}{\widetilde}

\newcommand{\bA}{\mathbf A}
\newcommand{\bx}{\mathbf x}
\newcommand{\bE}{\mathbf E}

\newcommand{\La}{\langle }
\newcommand{\Ra}{\rangle }
\newcommand{\rk}{\operatorname{rk}}

\newcommand{\tr}{\operatorname{tr}}

\newcommand{\e}{\varepsilon}
\newcommand{\be}{\begin{equation}}
\newcommand{\ee}{\end{equation}}

\newcommand{\fA}{\mathfrak{A}}

\newcommand{\ch}{\operatorname{ch}}

\newcommand{\nm }{\,\rule[-.6ex]{.13em}{2.3ex}\,}
\newcommand{\lnm}{\!\left\bracevert\!}
\newcommand{\rnm}{\!\right\bracevert\!}

\newcommand{\la}{\langle}
\newcommand{\ra}{\rangle}

\newcommand{\ti}[1]{_{\scriptstyle \text{\rm #1}}}
\newcommand{\ut}[1]{^{\scriptstyle \text{\rm #1}}}
\newcounter{vremennyj}

\numberwithin{equation}{section}

\newtheorem{thm}{Theorem}[section]
\newtheorem{lm}[thm]{Lemma}
\newtheorem{cor}[thm]{Corollary}

\newtheorem*{prop*}{Proposition}

\theoremstyle{remark}

\newtheorem*{rem*}{Remark}

\begin{document}

\title[Non-homogeneous vector square function]{Mixed $A_2$-$A_\infty$ estimates of the non-homogeneous vector square function with matrix weights}

\author{Sergei Treil}
 \thanks{Supported  in part by the National Science Foundation under the grant  DMS-1600139.}
\address{Department of Mathematics \\ Brown University \\ Providence, RI 02912 \\ USA}
\email{treil@math.brown.edu}



\begin{abstract}
This paper extends the results  from \cite{HyPeVo2017} about sharp $A_2$-$A_\infty$ estimates with matrix weights to the non-homogeneous situation.  
\end{abstract}

\maketitle

\section*{Notation}

\begin{entry}

\item[$\norm{\cdot}$, $\nm\,.\,\nm$]  norm; since we are dealing with vector- and
operator-valued functions, we will use the symbol $\|\,.\,\|$ (usually with a
subscript) for the norm in a function space, while $\nm\,.\,\nm$ is used for the
norm in the underlying vector (matrix) space, usually $\R^d$ (or $\C^d$) of the space of $d\times d$ matrices.
 Thus for a vector-valued function
$f$ the symbol
$\|f\|\ci{L^2}$ denotes its
$L^2$-norm, but the symbol $\nm f\nm$ stands for the scalar-valued function whose
value at a point $x$ is the norm of the vector $f(x)$; 

\item[$\la f\ra\ci I$] average, $\la f\ra\ci I = |I|^{-1} \int_I f (x) \dd x$; in agreement with the previous notation   for a vector-valued function $f$ the symbol $\la\nm f\nm\ra\ci I$ denotes the average of the function $x\mapsto \nm f(x)\nm$;
\item[$\bE\ci I$] averaging operator, $\bE\ci I f := \la f\ra\ci I \1\ci I$;
\item[$\ch I$] children of an interval $I$; 
\item[$\Delta\ci I$] Martingale difference operator, 
\[
\Delta\ci I := \sum_{I'\in\ch I} \bE\ci {I'} \quad -\qquad  \bE\ci I;
\]

\end{entry}

Expression $x\lesssim y$ means $x\le Cy$ with some absolute constant $C$; notation $x\underset{a,b}{\lesssim} y$ means that the constant $C$ depends only on parameters $a,b$

\section{Introduction and main results}

In this paper we obtain the result about sharp weighted estimates of vectors square function with matrix weights in the non-homogeneous  settings. In \cite{HyPeVo2017} such estimates were proved for the simplest case of homogeneous filtration, namely for the standard dyadic filtration on the real line $\R$. Formally,  our Corollary \ref{c: 1 weight sq funct} generalizes the main result from \cite{HyPeVo2017} from dyadic filtration on $\R$ to a general non-homogeneous filtration. 

We should mention that the result in \cite{HyPeVo2017} was proved only for the dyadic filtration on $\R$, and  its generalization even to the case of  dyadic lattice in $\R^N$ is not completely  trivial. The case of a non-homogeneous
filtration required some essential new ideas. In particular, a new construction of a dominating sparse square function, that takes care of non-homogeneity, was needed.   Moreover, \cite{HyPeVo2017} uses an estimate of the weighted sparse square function from  \cite{NPetTrV2017} that was proved  the reverse H\"{o}lder inequality for $A_\infty$ weights, which is not available in the non-homogeneous case.

Our main result (Theorem  \ref{t:Sq funct upper bd}) treat the two weight situation; while the corresponding two weigh result was not stated in \cite{HyPeVo2017}, the proof there works in the two weight case (but only for the dyadic filtration on $\R$).

\subsection{Setup}
\subsubsection{Atomic filtered spaces}
Let $(\mathcal{X}, \fA, \sigma)$ be a sigma-finite measure space with an atomic filtration $\fA_n$, that is, a sequence of increasing sigma-algebras $\fA_n\subset \cF$ such that for each $\fA_n$ there exists a countable collection $\mathcal{D}_n$ of disjoint sets of finite measure with the property that every set of $\fA_n$ is a union of sets in $\mathcal{D}_n$. 

We will often  use the notation $|A|$ for $\sigma(A)$ and $\dd x$ for $\dd \sigma(x)$. 

We will call the sets $I\in\cD_n$ \emph{atoms}, and denote by $\cD$ the collection of all atoms,  $\displaystyle\mathcal{D}=\cup_{n\in \mathbb{Z}}\mathcal{D}_n$. We allow a set $I$ to belong to several generations $\cD_n$, so formally an atom $I\in\cD_n$ is a pair $(I,n)$. To avoid overloading the notation, we skip the ``time'' $n$ and write $I$ instead of $(I,n)$; if we need to ``extract'' the time $n$, we will use the symbol $\rk I$. Namely, if $I$ denotes the atom $(I,n)$ then $n=\rk I$. 

The inclusion $I\subset J$ for atoms should be understood as inclusion for the sets together with the inequality $\rk I \ge \rk J$. However, the union (intersection) of atoms is just the union (intersection) of the corresponding sets and ``times'' $n$ are not taken into account. 

A standard example of such a filtration is the dyadic lattice $\mathcal{D}$ on $\mathbb{R}^N$, which explains the choice of notation. However, in what follows, $\mathcal{D}$ will always denote a general collection of atoms and $I\in \mathcal{D}$ will stand for an atom  in $\mathcal{D}$, and not necessarily for a dyadic interval. 

\subsubsection{Matrix-valued weights}
%Define
%\[
%\cF_0:=\{E\cap F: E\in \cF, \, F=\bigcup_1^n I_k, \ I_k\in\cD\}. 
%\] 
A matrix-valued weight $W$ is a locally integrable function on $\cX$ whose values are $d\times d$ positive semidefinite matrices. Here and everywhere in the paper \emph{locally integrable} means integrable on any atom $I$.   

The weighted space $L^2(W)$ is defined as the set of all measurable $\F^d$-valued functions ($\F=\R$ or $\C$) such that 
\[
\|f\|\ci{L^2(W)}^2 := \int_\cX \Bigl( W(x) f(x), f(x)\Bigr)_{\F^d} \dd x  <\infty ;
\]
as usual we take the quotient space over the set of functions of norm $0$.

\subsubsection{Matrix \texorpdfstring{$\bA_p$}{A_p} and \texorpdfstring{$\bA_\infty$}{A_infty} conditions} 
A matrix weight $W$ is said to satisfy the \emph{martingale} matrix $\mathbf A_2$ condition (and write $W\in(\mathbf A_2)$) if 
\[
[W]\ci{\bA_2} := \sup_{I\in\cD} \nm \la W\ra\ci I^{1/2} \la W^{-1}\ra\ci I^{1/2} \nm^2  < \infty\,.
\]
The quantity $[W]\ci{\bf A_2}$ is called the \emph{martingale $\mathbf A_2$ characteristic of the weight} $W$. 

Similarly for a pair of weights  $V$ and $W$ we   
the two-weight \emph{martingale} matrix $\bA_2$ condition 
\begin{align}
\label{2w-A_2}
\sup_{I\in\cD}  \nm\la W\ra\ci I^{1/2} \la V\ra\ci I^{1/2}\nm^2 =: [W, V]\ci{\bA_2} <\infty;
\end{align}

We say that a scalar weight $w$ on $\cX$ satisfies the \emph{martingale} $A_\infty$ condition if for all atoms $I\in\cD$ 
\begin{align}
\label{A_infty}
   \la M\ci I w  \ra\ci I \le C \la w\ra\ci I,
\end{align}
where $M\ci I$ is the localized to the atom $I$ maximal function
\begin{align}
\label{M_Q}
M\ci I f(x) = \sup \{ |\la f \ra\ci {I'} | \,: \, I'\in \cD(I), x\in \cX\}
\end{align}
(we put $M\ci I f(x) =  0$ for $x\notin I$). 

The best constant in \eqref{A_infty} is called the (martingale) $A_\infty$ characteristic of the weight $w$, and denoted by $[w]\ci{A_\infty}$. 

For a matrix weight $W$ define its \emph{scalar (martingale) $\bA_\infty$ characterisctic} $[W]\ci{\bA_\infty} = [W]\ci{\bA_\infty}\ut{sc}$ as 
\begin{align}
\label{A_infty-sc}
[W]\ci{\bA_\infty}%\ut{sc}
:= \sup_{e\in\F^d} [w_e ]\ci{A_\infty}, 
\end{align}
where the scalar weight $w_e$ is defined by $w_e(x) = (W(x) e, e)$, $x\in\cX $.

\subsubsection{Square functions}
For a weight $V$ the weighted square function $S^V$ is defined as 
\begin{align}
\label{eq: square function}
S^V f(x)  &= \biggl(\E_\alpha %\lnm \sum_{I\in\cD}  \alpha\ci I \Delta\ci I f (x) \rnm  \ 
 \biggl\bracevert\! V(x)^{1/2}\sum_{I\in\cD}  \alpha\ci I \Delta\ci I f (x) \!\biggr\bracevert^2 \biggr)^{1/2} \\
 \notag
& = 
\biggl( \sum_{I\in\cD}  \lnm V(x)^{1/2}  \Delta\ci I f (x) \rnm^2 \biggr)^{1/2}
\end{align}
The modified weighted square function $\wt S^V$ is given by 
\begin{align}
\label{eq: mod square function}
\wt S^V f(x) := \biggl(\sum_{I\in\cD} \sum_{I'\in\ch(I)} \lnm \la V\ra\ci {I'}^{1/2} \Delta\ci I f (x) \1\ci{I'}(x)\rnm^2 \biggr)^{1/2}\,. 
\end{align}
Integrating in $x$ we easily get 
\begin{lm}
\label{L: S vs mod S}
For any $f\in L^2(\cX, \F^d) = L^2(\F^d)$
\begin{align*}
\| S^V f \|^2\ci{L^2(\F^d)} = \| \wt S^V f \|^2\ci{L^2(\F^d)}
\end{align*}
\end{lm}

Recall that a collection $\cF\subset\cD$ is called $\e$-sparse if for any $I\in\cF$
\begin{align*}
\sum_{I'\in\ch\ci\cF(I) } |I'| \le \e|I|. 
\end{align*}
For a sparse family $\cF$ a weighted sparse \emph{averaging} square function $\cA^V=\cA^V\ci \cF$ is defined as 
\begin{align*}
\cA^V f (x) =\biggl( \sum_{I\in\cF} \left\la \lnm \la V\ra\ci I^{1/2} f \rnm  \right\ra\ci I^2 \1\ci I(x) \biggr)^{1/2}, \qquad f\in L^1\ti{loc}=L^1\ti{loc} (\cX;F^d) ;
\end{align*}
here, recall, $L^1\ti{loc}$ means integrable on any atom.

\subsection{Main results}
\begin{thm}
\label{t:Sq funct upper bd}
Let matrix weights $U$, $V$ satisfy the joint $\bA_2$ condition, and let $U\in(\bA_\infty)$. Then for any $f\in L^2=L^2(\F^d)$ 
\begin{align*}
\| S^V (U^{1/2} f) \|\ci{L^2} \underset{d}{\lesssim} [U,V]\ci{\bA_2}^{1/2} [U]\ci{\bA_\infty}^{1/2} \|f\|\ci{L^2}. 
\end{align*}
\end{thm}

In one weight situation this theorem gives us the following corollary.  

\begin{cor}
\label{c: 1 weight sq funct}
Let a matrix weight $W$ satisfy the $\bA_2$ condition. Then for any $f\in L^2(W)$
\begin{align*}
\| S^W f\|\ci{L^2} \underset{d}{\lesssim} [W]\ci{\bA_2}^{1/2} [W^{-1}]\ci{\bA_\infty}^{1/2} \|f\|\ci{L^2(W)}
\end{align*}
\end{cor}
To prove this corollary one just need to apply Theorem \ref{t:Sq funct upper bd} with $V=W$ and $U=W^{-1}$ to the function $W^{1/2} f$.

It is trivial that it is sufficient to prove Theorem \ref{t:Sq funct upper bd} only for the localized version $S^V\ci{I_0}$ of the square function $S^V$; in the localized version the sum in \eqref{eq: square function}   is taken only over $I\in\cD(I_0)$.  A uniform (in $I_0$ estimate for $S^V\ci{I_0}$ gives the same estimate for $S^V$

One can also consider the localized version $\wt S\ci{I_0}^V$ of $\wt S^V$, where again the sum in \eqref{eq: mod square function} is taken only over $I\in\cD(I_0)$. Clearly, Lemma \ref{L: S vs mod S} holds id we replace $S^V$ and $\wt S^V$ by $S^V\ci{I_0}$ and $\wt S^V\ci{I_0}$ respectively. 

The proof of Theorem \ref{t:Sq funct upper bd} is based on the following result about domination of a square function by a sparse one.

\begin{thm}
\label{t: Square sparse domination}
Let $\wt S^V =\wt S^V\ci{I_0}$ be the localized version of the square function \eqref{eq: mod square function}. 
Given $f\in L^1(\cX;\F^d)$ supported on $I_0$ a weight $V$
there exists $\frac12$-sparse family $\cF$ (depending on $f$ and $V$) such that 
\begin{align*}
\| \wt S^V f(x)\|\ci{L^2} \underset{d}{ \lesssim }  \| \cA^V f(x)\|\ci{L^2}\,.
\end{align*}
\end{thm}

For a measurable function $f$ we denote by $\La f\Ra\ci I$ its average, 
\[
\La f\Ra\ci I := \sigma(I)^{-1} \int_I f \dd\sigma;
\]
if $\sigma(I)=0$ we put $\La f\Ra\ci I=0$. The same definition is used for the vector and matrix-valued functions.  

In what follows we will often use $|E|$ for $\sigma(E)$ and $\dd x$ for $\dd\sigma$. 

To shorten the notation  we will denote by $\la\nm f\nm\ra\ci I$, where $f$ is a vector-valued function, the average of a function $x\mapsto \nm f(x)\nm$, 
\begin{align*}
\la\nm f\nm\ra\ci I := |I|^{-1} \int_I \nm f(x)\nm \dd x. 
\end{align*}
%%

%\begin{thm}
%\label{main}
%Let $W$ be a $d\times d$ matrix-valued weight and let $A\ci I$, $I\in\cD$ be a sequence of positive semidefinite $d\times d$ matrices. Then the following are equivalent:
%\begin{enumerate}
%\item $\displaystyle \sum_{I\in \mathcal{D}}\left\|A\ci I^{1/2}\La W^{1/2}f\Ra\ci I\right\|^2|I|\le A \|f\|^2\ci{L^2}$.
%\item $\displaystyle \sum_{I\in \mathcal{D}}\left\|A\ci I^{1/2}\La W f\Ra\ci I\right\|^2|I|\le A \|f\|^2_{L^2(W)}$.
%\item $ \displaystyle \frac{1}{|I_0|} \sum_{\substack{I\in\cD\\ I\subset I_0}}\La W\Ra \ci IA\ci I\La W\Ra \ci I |I|\le B\La W\Ra\ci{I_0} $ for all $I_0\in \mathcal{D}.$
%\end{enumerate}
%Moreover, $B\le A \le C B$, where $C=C(d)=e\cdot d^3(d+1)^2$.  
%%is a constant depending only on the dimension $d$. 
%\end{thm}
%

\section{Proof of the main results}

\subsection{Estimate of \texorpdfstring{$\cA^V$}{A<sup>V} } To prove Theorem \ref{t:Sq funct upper bd} we just need to estimate $\|\wt \cA^V(Uf)\|\ci{L^2}^2$. Clearly 
\begin{align*}
\|\cA^V(U^{1/2}f)\|\ci{L^2}^2 & = \sum_{I\in\cF} \Bigl\la \lnm \la V\ra\ci I^{1/2} U^{1/2}f \rnm  \Bigr\ra\ci I^2 |I| \\
&\le \sum_{I\in\cF} \lnm \la V\ra\ci I^{1/2}  \la U\ra\ci I^{1/2}\rnm^2 \Bigl\la \lnm \la U\ra\ci I^{-1/2} U^{1/2}f \rnm  \Bigr\ra\ci I^2 |I| 
\\
& \le [U,V]\ci{\bA_2} \sum_{I\in\cF}  \Bigl\la \lnm \la U\ra\ci I^{-1/2} U^{1/2}f \rnm  \Bigr\ra\ci I^2 |I| 
\end{align*}
So, we need to show
\begin{align*}
\sum_{I\in\cF}  \Bigl\la \lnm \la U\ra\ci I^{-1/2} U^{1/2}f \rnm  \Bigr\ra\ci I^2 |I| \underset{d}{\lesssim} [W]\ci{\bA_\infty}\ut{sc} \|f\|\ci{L^2}^2
\end{align*}
and for this it suffices to prove that for any scalar-valued $f\in L^2$
\begin{align}
\label{eq: needed bound for sparse}
\sum_{I\in\cF}  \Bigl\la \lnm \la U\ra\ci I^{-1/2} U^{1/2}\rnm |f|   \Bigr\ra\ci I^2 |I| 
\underset{d}{\lesssim }
[W]\ci{\bA_\infty} \|f\|\ci{L^2}^2 .
\end{align}

For $\alpha =\{\alpha\ci I\}\ci{I\in\cF} = \{|I|\}\ci{I\in\cF}$ consider the weighted space $\ell^2(\alpha) =\ell^2(\cF, \alpha)$, 
\[
\|\bx\|\ci{\ell^2(\alpha)}^2 = \sum_{I\in\cF} |x\ci I|^2 \alpha\ci I = \sum_{I\in\cF} |x\ci I|^2 |I|. 
\]
The estimate \eqref{eq: needed bound for sparse} is equivalent to the corresponding bound for the embedding operator $\cJ: L^2 \to \ell^2(\alpha)$ 
\[
\cJ f =\left\{\Bigl\la \lnm \la U\ra\ci I^{-1/2} U^{1/2}\rnm f  \Bigr\ra\ci I  \right\}_{I\in\cF} \,.
\]
The adjoint $\cJ^*:\ell^2(\alpha)\to L^2$ is given by 
\begin{align*}
\cJ^* \bx = \sum_{I\in\cF}  x\ci I \lnm \la U\ra\ci I^{-1/2} U^{1/2}\rnm \1\ci I,   
\end{align*}
meaning that 
\begin{align*}
(\cJ f, \bx)\ci{\ell^2(\alpha)} = (f, \cJ^* \bx)\ci{L^2}. 
\end{align*}
Since $\|\cJ^*\|\ci{\ell^2(\alpha)\to L^2}=\|\cJ\|\ci{L^2\to\ell^2(\alpha)}$, we reduced the problem to estimating the norm of $\cJ^*$. 

We can write 
\begin{align*}
\|\cJ^* \bx\|\ci{L^2}^2 = \sum_{I,J\in\cF} (\cJ\cJ^*)\ci{I,J} x\ci I x\ci J |I| |J|,
\end{align*}
where   $(\cJ\cJ^*)\ci{I,J} =0$ if $I\cap J=\varnothing$, $(\cJ\cJ^*)\ci{I,J} = (\cJ\cJ^*)\ci{J,I}\ge 0$ and for $I\subset J$
\begin{align*}
(\cJ\cJ^*)\ci{I,J} & =|J|^{-1}|I|^{-1}\int_I \lnm \la U\ra\ci I^{-1/2} U^{1/2}(x)\rnm \lnm \la U\ra\ci J^{-1/2} U^{1/2}(x)\rnm \dd x
\\
& \le |J|^{-1} \lnm \la U\ra\ci J^{-1/2}  \la U\ra\ci I^{1/2}\rnm  |I|^{-1}\int_I \lnm \la U\ra\ci I^{-1/2} U^{1/2}(x)\rnm^2 \dd x \\
& \underset{d}{\lesssim} |J|^{-1} \lnm  \la U\ra\ci J^{-1/2}  \la U\ra\ci I^{1/2}\rnm =: t\ci{I,J}
\end{align*}

Since $(\cJ\cJ^*)\ci{I,J} = (\cJ\cJ^*)\ci{J,I}\ge 0$
 it suffices to estimate $\sum_{I,J\in\cF} T\ci{I,J} |I| x\ci I |J|x\ci J$ where $x\ci I\ge 0$ and 
\begin{align*}
T\ci{I,J}= \left\{ \begin{array}{ll} t\ci{I,J} \qquad & I\subset J \\ 0 & I\not\subset J     \end{array}\right.
\end{align*}
We will need the following well-known result, see \cite[Lecture VII]{Nik_shift_1986} 

\begin{lm}[Senichkin--Vinogradov test, AKA iterated kernel test]
\label{l:Vin-Sen}
Let $k(\fdot, \fdot)\ge 0$ be a measurable locally integrable%
\footnote{Here ``locally integrable means that the right hand side of \eqref{l:Vin-Sen} is finite for some dense collection of non-negative functions $f$. In our case the measure is a discrete one, so any kernel is locally integrable.}
 function on $X\times X$. Let 
\begin{align*}
\int_X k(s,x) k(s,t) \dd\mu(s) \le C [k(x,t) + k(t,x)] . 
\end{align*}
Then for any measurable $f\ge 0$
\begin{align}
\label{eq: V-S result}
\iint_{X\times X} k(s,t) f(s) f(t) \dd\mu(s)\dd\mu(s) \le 2C\|f\|\ci{L^2(\mu)}^2
\end{align}
\end{lm}

Let us apply the above Lemma \ref{l:Vin-Sen}. We need to show that for $J,K\in\cF$, $J\subset K$
\begin{align}
\label{eq:check Vin-Sen}
\sum_{I\in\cF(J)} t\ci{I,J} t\ci{I,K} |I| \underset{d}{\lesssim} [U]\ci{\bA_\infty} t\ci{J,K}.
\end{align}
A simple calculation gives us 
\begin{align*}
\sum_{I\in\cF(J)} t\ci{I,J} t\ci{I,K} |I| &= |K|^{-1}|J|^{-1} \sum_{I\in\cF(J)}  \lnm  \la U\ra\ci J^{-1/2}  \la U\ra\ci I^{1/2}\rnm  \lnm  \la U\ra\ci K^{-1/2}  \la U\ra\ci I^{1/2}\rnm |I|
\\ &
\le |K|^{-1}|J|^{-1} \lnm  \la U\ra\ci K^{-1/2}  \la U\ra\ci J^{1/2}\rnm \sum_{I\in\cF(J)} 
\lnm  \la U\ra\ci J^{-1/2}  \la U\ra\ci I^{1/2}\rnm^2 |I| .
\end{align*}
Continuing the estimates we write
\begin{align*}
\sum_{I\in\cF(J)} 
\lnm  \la U\ra\ci J^{-1/2}  \la U\ra\ci I^{1/2}\rnm^2 |I| 
& \le \sum_{I\in\cF(J)} 
\tr\left\{  \la U\ra\ci J^{-1/2}  \la U\ra\ci I \la U\ra\ci J^{-1/2} \right\} |I|
\\
& = \sum_{I\in\cF(J)} \tr \la \wt U\ra\ci I  |I|
= \sum_{I\in\cF(J)}  \la \tr \wt U\ra\ci I  |I|, 
\end{align*}
where the weight $\wt U$ is given by $ \la U\ra\ci J^{-1/2} U  \la U\ra\ci J^{-1/2}$, so $\la \wt U\ra\ci J =\bI$. 

The scaling $U\mapsto A^* U A$ with a constant invertible matrix $A$ does not change the  $\bA_\infty$ constant of the weight. Therefore the scalar weight $u:=\tr \wt U$ is a scalar $A_\infty$ weight with $[u]\ci{A_\infty} \le [\wt U ]\ci{\bA_\infty} = [ U ]\ci{\bA_\infty}$, as the sum of $d$ scalar $A_\infty $ weighs, see \cite[Lemma 4.5]{NPetTrV2017}. 

It is an easy and well-known fact that for an $\e$-sparse family
\begin{align*}
|J|^{-1}\sum_{I\in\cF(J)} \le (1-\e)^{-1} \la M\ci J u\ra\ci J, 
\end{align*}
so 
\begin{align*}
\sum_{I\in\cF(J)} 
\lnm  \la U\ra\ci J^{-1/2}  \la U\ra\ci I^{1/2}\rnm^2 |I| & \le 2 \la M\ci J u\ra\ci J|J|
\\
& \le 2 [U]\ci{\bA_\infty} \la u \ra\ci J |J|
\\
& = 2 [U]\ci{\bA_\infty} (\tr\bI ) |J| = 2d [U]\ci{\bA_\infty} |J|.
\end{align*}
Gathering everything together we get \eqref{eq:check Vin-Sen} with the implied constant $2d$. Thus the estimate \eqref{eq: needed bound for sparse}, and so Theorem  \ref{t:Sq funct upper bd}, are proved. \hfill\qed

\subsection{Proof of sparse domination (Theorem \ref{t: Square sparse domination}) }

The proof is pretty standard, just with a few twists. 

Before starting, let us make a trivial observation, that the problem is invariant to the scaling by a constant invertible matrix. Namely, if $A$ is a constant invertible matrix, $\wt V:= A^* V A$, $\tilde f := A^{-1} f$, then 
\begin{align}
\label{eq: rescaling S}
\wt S^V f(x) = \wt S^{\wt V} \tilde f (x), \qquad \cA^{V} f (x) = \cA^{\wt V} \tilde f (x). 
\end{align}
So, to simplify the notation we will on each step of the induction construction of the sparse family $\cF$ rescale the weight $V$ and the function $f$. This rescaling does not change anything, it just makes the formulas for stopping moments simpler; they of course can be equivalently rewritten without rescaling, but the formulas will be uglier.  

We start with the atom $I_0$, which will be initial atom  in the family $\cF$. Let us rescale $V$ anf $f$, defining
\begin{align*}
\tilde f := \la V \ra\ci{I_0}^{1/2} f, \qquad \wt V := \la V \ra\ci{I_0}^{-1/2} V \la V \ra\ci{I_0}^{-1/2}, 
\end{align*}
so $\la \wt V\ra\ci{I_0}=\bI$. 

We then pick the collection of stopping intervals $\cF^1(I_0)$ which consists of \emph{maximal} by inclusion intervals $J\in\cD(I_0)$ for which either one of the following 3 conditions holds
\begin{align*}
\sum_{I\in\cD(I_0): I \supsetneqq J} \lnm \Delta\ci I \tilde f\rnm^2 > C^2 \la \nm \tilde f \nm \ra\ci{I_0}^2, 
\qquad \tr \la \wt V \ra\ci I  > C d, \qquad \la \nm\tilde f \nm \ra\ci I > C\la \nm\tilde f \nm \ra\ci {I_0}
\end{align*}
Note that by picking sufficiently large $C=C_k(\e)$ we can assure that the total measure of maximal intervals where one of the condition holds is at most $\e|I_0|$. For the third and second condition it follows from the weak type estimates for the maximal function, and for the first condition from the weak type estimates for the square function.%
\footnote{Note that the weak type estimates hold for martingales with values in a Hilbert space, see \cite[Theorem 3.6]{Burkholder_LNM-1991} , so the constant $C$ does not depend on $d$. Of course, for the finite-dimensional case one can easily get the weak type estimates from the scalar result (with the constant depending on $d$).  }
Thus picking $\e=1/6$ we get that 
\begin{align*}
\sum_{J\in \cF^1(I_0) }|J| \le |I_0|/2. 
\end{align*}

It is trivial that for $x\in I_0\setminus  \bigcup_{J\in\cF^1(I_0)} J$ we have 
\begin{align*}
\left(\wt S^{\wt V} \tilde f (x) \right)^2\le C^3 d \la \nm \tilde f \nm \ra\ci{I_0}^2. 
\end{align*}

Now take $J\in\cF^1(I_0)$ and for $I\in\cD$ let $\hat I$ be  its parent. Then for $x\in J$
\begin{align*}
\left(\wt S^{\wt V}\tilde f (x)\right)^2 = \sum_{\substack{I\in\cD(I_0): x\in I\\ J\subsetneqq I\subsetneqq I_0} }  \lnm \la \wt V \ra\ci I^{1/2} (\bE\ci I - \bE\ci{\hat I})\tilde f \rnm^2
+ \lnm \la \wt V \ra\ci J^{1/2} (\bE\ci J - \bE\ci{\hat J})\tilde f \rnm^2
+ \left( \wt S^{\wt V}\ci{J} \tilde f (x) \right)^2
\end{align*}
The first term is clearly estimated by $ C^3 d \la \nm \tilde f \nm \ra\ci{I_0}^2$. To estimate the second term we write
\begin{align*}
\lnm \la \wt V \ra\ci J^{1/2} (\bE\ci J - \bE\ci{\hat J})\tilde f \rnm  &\le
\la \lnm \la\wt V \ra\ci J^{1/2} \tilde f \rnm\ra\ci J + \la \nm \tilde f\nm \ra\ci{\hat J} \lnm \la \wt V \ra\ci J^{1/2}\rnm 
\\
& \le \la \lnm \la\wt V \ra\ci J^{1/2} \tilde f \rnm\ra\ci J + C \la \nm \tilde f\nm \ra\ci{I_0} 
\lnm \la \wt V \ra\ci J^{1/2}\rnm.
\end{align*}
Gathering everything together we get that
\begin{align*}
\left(\wt S^{\wt V}\tilde f (x)\right)^2 
&\le 2C^3 d 
\biggl(     \la \nm \tilde f \nm \ra\ci{I_0}^2  
 + 
\sum_{J\in\cF^1(I_0)} \la \lnm \la\wt V \ra\ci J^{1/2} \tilde f \rnm\ra\ci J^2 \1\ci J
\\
& \qquad \qquad  +     
\quad \sum_{J\in\cF^1(I_0)} \la \nm \tilde f\nm \ra\ci{I_0}^2 
\lnm \la \wt V \ra\ci J^{1/2}\rnm^2 \1\ci J 
\biggr)
+  \quad \sum_{J\in\cF^1(I_0)} \left( \wt S\ci J^{\wt V} \tilde f(x)\right)^2
\\
& = 2C^3 d 
\biggl(     \la \nm \la V \ra\ci {I_0}^{1/2} f \nm \ra\ci{I_0}^2  
 + 
\sum_{J\in\cF^1(I_0)} \la \lnm \la V \ra\ci J^{1/2}  f \rnm\ra\ci J^2 \1\ci J
\\ &
 \qquad \qquad   +     
\quad \sum_{J\in\cF^1(I_0)} 
\la \nm \la V \ra\ci {I_0}^{1/2} f \nm \ra\ci{I_0}^2 
\nm \la V\ra\ci{I_0}^{-1/2}\la  V \ra\ci J^{1/2}\nm^2 \1\ci J 
\biggr)
+  \quad \sum_{J\in\cF^1(I_0)} \left( \wt S\ci J^{ V}  f(x)\right)^2
\end{align*}

Repeating this procedure with each $\wt S^V\ci J f $, $J\in \cF^1(I_0)$ and iterating we get that 
\begin{align*}
\wt S^{\wt V}\tilde f (x) \le C d^{1/2} \left( \cA^V f (x) + \cA^V\ti{mod} f (x) \right) 
\end{align*}
(here we do not write $\cF$ as the index, but keep in mind that tho operators depend on $\cF$), where 
\begin{align}
\label{eq: A_mod}
\left( \cA^V\ti{mod} f (x) \right)^2 = 
\sum_{I\in\cF} \sum_{I'\in\ch(I)} \la \nm \la V \ra\ci {I}^{1/2} f \nm \ra\ci{I}^2 
\nm \la V\ra\ci{I}^{-1/2}\la  V \ra\ci{I'}^{1/2}\nm^2 \1\ci{I'} 
=: \sum_{I\in\cF} \sum_{I'\in\ch(I)} F\ci I.
\end{align}
If we estimate the norm of each term in \eqref{eq: A_mod} we get that
\begin{align*}
\|F\ci I\|\ci {L^1}
& \le 
\la \nm \la V \ra\ci {I}^{1/2} f \nm \ra\ci{I}^2 \biggl\|\sum_{I'\in\ch(I)} 
\nm \la V\ra\ci{I}^{-1/2}\la  V \ra\ci{I'}^{1/2}\nm^2 \1\ci{I'}\biggr\|_{L^1}
\\
& \le \la \nm \la V \ra\ci {I}^{1/2} f \nm \ra\ci{I}^2 
\sum_{I'\in\ch(I)} \tr\left( \la V\ra\ci{I}^{-1/2}\la  V \ra\ci{I'} \la V\ra\ci{I}^{-1/2}\right) |I'|
\\
& \le \la \nm \la V \ra\ci {I}^{1/2} f \nm \ra\ci{I}^2 |I| d, 
\end{align*}
so it is dominated by the norm of the corresponding term in $\cA^V$.
The theorem is proved. \hfill\qed

\def\cprime{$'$}
  \def\lfhook#1{\setbox0=\hbox{#1}{\ooalign{\hidewidth\lower1.5ex\hbox{'}\hidewidth\crcr\unhbox0}}}
\providecommand{\bysame}{\leavevmode\hbox to3em{\hrulefill}\thinspace}
\providecommand{\MR}{\relax\ifhmode\unskip\space\fi MR }
\providecommand{\MRhref}[2]{%
  \href{http://www.ams.org/mathscinet-getitem?mr=#1}{#2}
}
\providecommand{\href}[2]{#2}

\end{document}